\documentclass{article}
\usepackage{graphicx}
\usepackage{amsmath}
\usepackage{amsfonts}
\usepackage[indent=10pt]{parskip}
\usepackage[english]{babel}
\usepackage{babel}
\usepackage{csquotes}
\usepackage{hyperref}
\usepackage{amsthm}
\usepackage{inputenc}
\usepackage{biblatex}
\newtheorem{theorem}{Theorem}
\newtheorem{corollary}[theorem]{Corollary}
\newtheorem{proposition}[theorem]{Proposition}
\newtheorem{lemma}[theorem]{Lemma}
\newtheorem*{remark}{Remark}
\usepackage{biblatex} 
\DeclareMathOperator{\esssup}{esssup}

\title{On equality of the $L^\infty$ norm of the gradient under the Hausdorff and Lebesgue measure}
\author{Ng Ze-An}
\date{}

\begin{document}
\maketitle
\begin{abstract}
    Let $\Omega$ be an open subset of $\mathbb R^n$, and let $f: \Omega \to \mathbb R$ be differentiable almost everywhere with respect to $s$-dimensional Hausdorff measure $\mathcal H^s$, where $0 \leq s \leq n-1$. We show that $\|\nabla f\|_{L^\infty (\mathcal H^s)} = \|\nabla f\|_{L^\infty(\mathcal H^n)}.$ We deduce that convergence in the Sobolev space $W^{1, \infty}$ preserves everywhere differentiability.
\end{abstract}

\section{Introduction}
Differentiability almost everywhere is a cornerstone of modern analysis. The celebrated Rademacher's theorem asserts that every Lipschitz function on $\mathbb R^n$ is differentiable almost everywhere with respect to Lebesgue measure. In many settings, however, particularly in geometric measure theory, partial differential equations and the calculus of variations, one encounters functions that are differentiable almost everywhere with respect to finer, lower dimensional measures. 

For instance, Evans and Smart proved in \cite{evanssmart} that solutions to the $\infty$-Laplace equation are everywhere differentiable, which may be viewed as almost everywhere differentiability with respect to $0$-dimensional Hausdorff measure. For more general functions, differentiability may be known everywhere except for a thin set, such as a manifold or a rectifiable set, or even a fractal set of non-integer Hausdorff dimension.

In such situations, a natural question arises - how much can one trust the $L^\infty$ norm of the gradient computed with respect to the classical Lebesgue measure? Suppose a continuous function is differentiable almost everywhere with respect to $s$-dimensional Hausdorff measure for some $s < n$. Since the $L^\infty$ norm with respect to Lebesgue measure ignores all Lebesgue null sets, which may still be large with respect to the Hausdorff measure, it is conceivable that the largest values of $|\nabla f|$ may hide on a set invisible to Lebesgue measure, but visible to Hausdorff measure. 

Our main theorem says that for $0 \leq s \leq n-1$, this is not the case for functions that are differentiable $\mathcal H^s$ almost everywhere. That is, the $L^\infty$ norm of the gradient is the same under the Hausdorff and Lebesgue measure.

\begin{theorem}\label{maintheorem}
Let $\Omega$ be an open subset of $\mathbb R^n$, and let $f: \Omega \to \mathbb R$ be continuous, and further differentiable $\mathcal H^s$-almost everywhere, for some real $0 \leq s \leq n-1$. Then we have

\begin{equation}\label{eq1}\sup_{\Omega} |\nabla f| = \underset{\Omega}{\esssup} \, |\nabla f|\end{equation}

\noindent where the essential supremum is taken with respect to Lebesgue measure. In particular,
\begin{equation}\label{eq2}\|\nabla f\|_{L^\infty (\mathcal H^s)} = \|\nabla f\|_{L^\infty(\mathcal H^n)}.\end{equation}
\end{theorem}
\begin{remark}Note that in the statement of Theorem \ref{maintheorem}, $\mathcal H^n$ corresponds to the usual Lebesgue measure on $\Omega$.\end{remark}

 The hypothesis that the exceptional set of non-differentiability be of strictly smaller dimension is crucial. Under only the hypothesis of differentiability almost everywhere with respect to Lebesgue measure, the equality \ref{eq1} need not hold, as demonstrated in \cite{FERNANDEZSANCHEZ2016713}. In this work, the authors construct a singular function on $\mathbb R$ that is differentiable with derivative $0$ Lebesgue almost everywhere, yet takes a nonzero derivative on a set of Hausdorff dimension $1$. In particular, the $L^\infty$ norm of the gradient under the Lebesgue measure and any lower dimensional Hausdorff measure differs.

We deduce from Theorem \ref{maintheorem} some functional analytic consequences. The first concerns preservation of everywhere differentiability under $W^{1, \infty}$ convergence. For $1 \leq p < \infty$, convergence in $W^{1, p}$ is too weak to preserve pointwise classical differentiability.

In contrast, we show as a corollary of Theorem \ref{maintheorem} that if a sequence of everywhere differentiable functions converges in $W^{1, \infty}$, then the limit is again everywhere differentiable, with uniform convergence of the gradients.

\begin{corollary}\label{corollary2}
Let $f_n: \Omega \to \mathbb R$ be everywhere differentiable. Assume that $f_n \to f$ in $W^{1, \infty}$. Then $f$ is everywhere differentiable and further $f_n \to f$ and $\nabla f_n \to \nabla f$ uniformly.
\end{corollary}

As mentioned in the opening paragraphs, in the $L^\infty$ calculus of variations (see \cite{viscosity} or \cite{lindqvist2016notes} for an overview), one often encounters $W^{1, \infty}$ functions that are differentiable in some stronger sense than Lebesgue almost everywhere. One is also often faced with $W^{1, \infty}$ limits of such functions. Thus Corollary \ref{corollary2} may be of intrinsic interest in this field.

In preparation for our next corollary, recall that if $1 \leq p < \infty$, the closure of $C^1 (\Omega)$ in $W^{1, p}$ is the whole space. In sharp contrast, our next corollary shows that $C^1 (\Omega)$ is closed in $W^{1, \infty}$, revealing another qualitative difference between the $L^p$ and $L^\infty$ setting.

\begin{corollary}\label{corollary3}
$C^1 (\Omega)$ is a closed subspace of $W^{1, \infty} (\Omega)$.
\end{corollary}
We now describe how the rest of the paper is structured. In Section 2, we present the proof of the main Theorem \ref{maintheorem}. In Section 3, we present the proof of Corollary \ref{corollary2} and \ref{corollary3}. In Section 4, we present some natural questions and directions for further research.

\section{Proof of main theorem}
We begin by proving the case with $s = 0$ and $n = 1$, restated below for the readers' convenience.

\begin{proposition}\label{oned}
   Let $\Omega$ be an open subset of $\mathbb R$, and let $f: \Omega \to \mathbb R$ be differentiable everywhere. Then 
   $$\underset{\Omega}{\sup} \, |f'| = \underset{\Omega}{\esssup} \, |f'|.$$
\end{proposition}
\begin{proof}
As every open subset in $\mathbb R$ is a countable union of open intervals, without loss of generality, we may work on an open interval $I$ instead. 

 Since $$\sup_I |f'|=\sup_{x, y \in I, x < y} \Big|\frac{f(y)-f(x)}{y-x}\Big|,$$ it will suffice to show the following mean value theorem for everywhere differentiable functions - for all intervals $[a, b] \subset I$,
\begin{equation}\label{meanvalue} 
 \frac{f(b)-f(a)}{b-a}\leq \esssup_{[a, b]} f' .
\end{equation}
To this end, write $v(x):=f(x)-\frac{f(b)-f(a)}{b-a}(x-a)$. Then (\ref{meanvalue}) is equivalent to the claim that $\esssup_{[a, b]} v' \geq 0$. 

\noindent Assume for contradiction the latter did not hold. Then we would have that $v$ is everywhere differentiable with $v'(x) < 0$ almost everywhere. We claim this implies that $v$ is (non-strictly) decreasing on $(a, b)$. Since $v(a) = v(b) = f(a)$, this would imply that $v$ is necessarily constant, which is the desired contradiction.

To see that $v$ is necessarily decreasing, we again proceed by contradiction. So assume $v$ is non-decreasing, so $v(z)>v(y)$ for some $z > y$.

The intervals $[c,d]\subset(y, z)$ with $v(d)<v(c)$ cover the set $\{v'<0\} \cap (y,z)$ in the Vitali sense, that is,  every point of $\{v'<0\} \cap (y,z)$ belongs to arbitrarily small such intervals. By the Vitali covering lemma, there is a finite  disjoint family of these intervals whose sum of lengths is greater than $\frac12(z-y)$. In other words, labelling these intervals $[c_{2k-1},c_{2k}]$ for $k=1\dots n$, there exists  a finite  sequence $$y =: c_0 < c_1<\dots < c_{2n+1}:=z$$  such that $$v(c_{2k})<v(c_{2k-1})$$ and

$$
\begin{aligned}
\sum_{k=0}^n(c_{2k+1}-c_{2k})=(z-y) &- \sum_{k=1}^n(c_{2k}-c_{2k-1})  \\
&\le \frac12(z-y).
\end{aligned}$$
\noindent Therefore 
$$\begin{aligned}
v(z)-v(y)&= \sum_{j=0}^{2n} v(c_{j+1})-v(c_j)  \\
&\le \sum_{k=0}^nv(c_{2k+1})-v(c_{2k}) \\
&=\sum_{k=0}^n\frac{v(c_{2k+1})-v(c_{2k})}{c_{2k+1}-c_{2k}}(c_{2k+1}-c_{2k}) \\
&\le  \sum_{k=0}^n(c_{2k+1}-c_{2k})\max_{0\le k\le n} \frac{v(c_{2k+1})-v(c_{2k})}{c_{2k+1}-c_{2k}} \\
&\le\frac12(z-y)\max_{0\le k\le n} \frac{v(c_{2k+1})-v(c_{2k})}{c_{2k+1}-c_{2k}} \\
\end{aligned}$$

\noindent That is, for the maximizing index $k^*$, the interval $[c,d]:=[c_{2k^*+1},c_{2k^*}]$ has  $$\frac{v(d)-v(c)}{d-c} \ge2 \frac{v(z)-v(y)}{z-y}.$$ If we iterate this procedure we get a nested sequence $[a_n,b_n]\subset [a,b]$ with $\frac{v(b_n)-v(a_n)}{b_n-a_n}\to\infty.$ If $x_*\in\bigcap_{n\ge0}[a_n,b_n]$, we have $\limsup_{x\to x^*}\frac{v(x^*)-v(x)}{x^*-x}=+\infty$, a contradiction to the fact that $f$, hence $v$ is everywhere differentiable, in particular at $x^*$, so $\limsup_{x\to x^*}\frac{v(x^*)-v(x)}{x^*-x}=v'(x^*) < +\infty$.
\end{proof}

Now we turn to the proof of the general case. Before we do so, we require some preliminaries on spherical coordinates. 

We let $\Phi: (r, \theta_1, \dots, \theta_{n-1}) \to \mathbb R^n$ be standard spherical coordinates on $\mathbb R^n$. Explicitly,

$$\begin{aligned}
\Phi_1 &= r \cos\theta_1, \\[4pt]
\Phi_2 &= r \sin\theta_1 \cos\theta_2, \\[4pt]
\Phi_3 &= r \sin\theta_1 \sin\theta_2 \cos\theta_3, \\[4pt]
&\ \,\vdots \\[4pt]
\Phi_{n-1} &= r \sin\theta_1 \sin\theta_2 \cdots \sin\theta_{n-2}\cos\theta_{n-1}, \\[4pt]
\Phi_n &= r \sin\theta_1 \sin\theta_2 \cdots \sin\theta_{n-2}\sin\theta_{n-1}.
\end{aligned}$$
The variables range as follows - $r \in (0, \infty)$, $\theta_1, \dots, \theta_{n-2} \in [0, \pi]$ and $\theta_{n-1} \in [0, 2\pi]$.

We note the following general fact about spherical coordinates, which states that the restriction of the coordinate system to a slice with the last angular coordinate $\theta_{n-1}$ fixed yields a hyperplane through the origin.

\begin{lemma}\label{hypersurface}
Fix $z \in [0,2\pi)$, and let $S_z := \Phi((0, \infty] \times[0, \pi]^{n-2} \times \{z\})$ denote the subset of $\mathbb{R}^n$
obtained by holding the final angular coordinate $\theta_{n-1}$ equal to $z$
in the standard spherical coordinate parametrization.
Then $S_z$ coincides with the hyperplane
\[
H_z = \bigl\{\theta \in \mathbb{R}^n :\,
\theta_{n-1}\sin z - \theta_n\cos z = 0 \bigr\}.
\]
In particular, $H_z$ is a linear subspace of codimension one passing through the origin.
\end{lemma}
\begin{proof}
When the final angular coordinate is fixed at $\theta_{n-1}=z$,
the last two components of the spherical coordinate map take the form
\[
(\Phi_{n-1},\Phi_n)
   = r\,P(\theta_1,\ldots,\theta_{n-2})
      \bigl(\cos z,\sin z\bigr),
\]
where $P(\theta_1,\ldots,\theta_{n-2})
      = \sin\theta_1\sin\theta_2\cdots\sin\theta_{n-2}$.
Thus every point with $\theta_{n-1}=z$ satisfies the linear relation
\[
x_{n-1}\sin z - x_n\cos z = 0.
\]
The set of all such points is precisely the hyperplane $H_z$,
which is homogeneous and therefore passes through the origin. This $S_z \subset H_z$.

Conversely, let $x \in H_z$.  The defining relation
$x_{n-1}\sin z - x_n\cos z = 0$ implies that
$(x_{n-1},x_n) = \lambda(\cos z,\sin z)$ for some $\lambda \ge 0$.
Setting $r = |x|$ and $P = \lambda/r \in [0,1]$, the remaining coordinates
satisfy $x_1^2 + \cdots + x_{n-2}^2 = r^2(1-P^2)$.

To recover the corresponding angular parameters, define them recursively as follows:
\[
\theta_1 = \arccos\!\left(\frac{x_1}{r}\right),
\quad
\theta_2 = \arccos\!\left(
   \frac{x_2}{r\sin\theta_1}
\right),
\quad
\theta_3 = \arccos\!\left(
   \frac{x_3}{r\sin\theta_1\sin\theta_2}
\right),
\]
and in general,
\[
\theta_j = \arccos\!\left(
   \frac{x_j}{r\sin\theta_1\cdots\sin\theta_{j-1}}
\right),
\quad
1 \le j \le n-2.
\]
By construction, these angles lie in $(0,\pi)$ and satisfy
\[
r\sin\theta_1\cdots\sin\theta_{j-1}\cos\theta_j = x_j
\quad (1 \le j \le n-2),
\qquad
r\sin\theta_1\cdots\sin\theta_{n-2} = \lambda.
\]
With these choices and $\theta_{n-1}=z$, we have
\[
\Phi(r,\theta_1,\ldots,\theta_{n-2},z)
   = (x_1,\ldots,x_{n-2},\lambda\cos z,\lambda\sin z)
   = x,
\]
showing that every $x\in H_z$ arises as
$\Phi(r,\theta_1,\ldots,\theta_{n-2},z)$.
Hence $H_z \subseteq S_z$.

We will need two more lemmas before the proof of Theorem \ref{maintheorem}. The first is a slicing lemma for Hausdorff measure.

\begin{lemma}\label{slicing}Let $1 \leq k \leq n$ be an integer, and let $E \subset \mathbb R^n$ be an arbitrary set of null $\mathcal H^k$ measure. Then for $\mathcal H^1$-almost every  $z \in [0, 2\pi]$, we have $\mathcal H^{k-1} (E \cap S_{z}) = 0$.\end{lemma} 

\begin{proof} Since $\Phi$ is a locally bi-Lipschitz map, we have $\mathcal H^u(E) = 0$ iff $\mathcal H^u (\Phi^{-1}(E)) = 0$, for any $0 \leq u \leq n$ and measurable set $E$. 

\noindent Setting $J_{z} = \Phi^{-1} (S_{z})$, we claim that
$$\mathcal H^{k-1} (J_z \cap \Phi^{-1} (E)) = 0.$$
for $\mathcal H^{1}$ almost every $z \in [0, 2\pi]$.

\noindent Indeed, $J_z$ is simply the rectangular slice $(0, \infty) \times [0, \pi]^{n-2} \times \{z\}$. Writing the condition $\mathcal H^k(\Phi^{-1} (E)) = 0$ in integral form, and expanding using the coarea formula for Hausdorff measures, see \cite{coarea}, Theorem 1.23, we have
$$\begin{aligned} 0 &= \int_{(0, \infty] \times [0, \pi]^{n-2} \times [0, 2\pi]} \mathbf 1_{\Phi^{-1}(E)} \, d \mathcal H^k \\
&= \int_{[0, 2\pi]} \int_{(0, \infty] \times [0, \pi]^{n-2}} \mathbf 1_{\Phi^{-1}(E)}(w,z) \, d\mathcal H^{k-1}(w) d\mathcal H^1 (z), \end{aligned}$$
\noindent and since the integrand is nonnegative, the integral being zero implies that the inner integral is $0$ for $\mathcal H^1$ a.e. $z \in [0, 2\pi]$. But the inner integral is exactly $\mathcal H^{k-1} (J_z \cap \Phi^{-1}(E))$. Since $E \cap S_z = \Phi(J_z \cap \Phi^{-1} (E))$, the conclusion of the lemma follows from the statement at the beginning of the proof.\end{proof}

\noindent The second lemma is a Fubini type equality for iterated $L^\infty$ norms.

\begin{lemma}\label{iterated}
\noindent With $S_z$ as above, we have for any $g \in L^\infty (\mathcal H_n, B_R)$,
$$\underset{ \mathcal H^1, z \in [0, 2\pi]}{\esssup} \|\nabla f\|_{L^\infty (\mathcal H^{n-1}, S_z \cap B_R(0))} = \|\nabla f\|_{L^\infty (\mathcal H^n, B_R(0))}.$$
\end{lemma}

\begin{proof}
\noindent We define the exceptional set
$$B := \{|\nabla f| > \underset{ \mathcal H^1, z \in [0, 2\pi]}{\esssup} \|\nabla f\|_{L^\infty (\mathcal H^{n-1}, S_z \cap B_R(0))}\} \cap B_R(0).$$
\noindent We claim that $\mathcal H^n (B) = 0$. By definition of the $L^\infty$ norm, this would imply 
$$\|\nabla f\|_{L^\infty (\mathcal H^n, B_R(0))} \leq \underset{ \mathcal H^1, z \in [0, 2\pi]}{\text{esssup}} \|\nabla f\|_{L^\infty (\mathcal H^{n-1}, S_z \cap B_R(0))}.$$
\noindent Applying the bi-Lipschitz property of $\Phi$ as in the proof of Lemma 6, it is equivalent to show that $\mathcal H^{n} (\Phi^{-1}(B)) = 0$. But we may write
$$\begin{aligned}\mathcal H^n (\Phi^{-1}(B)) & = \int_{[0, R] \times [0, \pi]^{n-2} \times [0, 2\pi]} \mathbf 1_{\Phi^{-1}(B)} \, d\mathcal H^n \\
&= \int_{[0, 2\pi]} \int_{J_z} \mathbf 1_{\Phi^{-1}(B)} (w, z) \, d\mathcal H^{n-1} (w) \, d\mathcal H^1 (z) \\ 
& = \int_{[0, 2\pi]}\mathcal H^{n-1}(\Phi^{-1}(B) \cap J_z \cap B_R(0)) \, d\mathcal H^1.\end{aligned}$$

\noindent By assumption, we have $\mathcal H^{n-1}(B \cap S_z) = 0$ for $\mathcal H^1$ a.e. $z \in [0, 2\pi]$. hence $\Phi^{-1}(B) \cap J_z = \Phi^{-1} (B \cap S_z)$ has $\mathcal H^{n-1}$ measure $0$ as well for $\mathcal H^1$ a.e. $z$. Thus
$$\int_{[0, 2\pi]}\mathcal H^{n-1}(\Phi^{-1}(B) \cap J_z) \, d\mathcal H^1 = 0,$$
\noindent and we conclude.

\noindent The reverse inequality follows similar lines. Write 

$$D := \{|\nabla f| > \|\nabla f\|_{L^\infty (\mathcal H^n, B_R(0))}\}.$$

\noindent The desired inequality then follows if we can show $\mathcal H^{n-1} (S_z \cap D) = 0$ for $\mathcal H^1$ a.e. $z$.

\noindent But we have
$$\begin{aligned}0 &= \mathcal H^n (D) \\
&= \mathcal H^n (\Phi^{-1}(D))\\
&= \int_{[0, 2\pi]} \mathcal H^{n-1}(J_z \cap \Phi^{-1}(D)) \, d\mathcal H^1 (z),\end{aligned}$$
\noindent which implies that  $\mathcal H^{n-1}(J_z \cap \Phi^{-1}(D)) = 0$ for $\mathcal H^1$ a.e. $z \in [0, 2\pi]$. Since $S_z \cap D = \Phi( J_z \cap \Phi^{-1}(D))$, we conclude by the Lipschitz property of $\Phi$.
\end{proof}

We now present the proof of Theorem \ref{maintheorem}.

\end{proof}
\begin{proof}[Proof of Theorem 1]

Let $x \in \Omega$ be an arbitrary point of differentiability of $f$. We show that $|\nabla f(x)| \leq \|\nabla f\|_{L^\infty (\mathcal H^n, \Omega)}$. This will imply (\ref{eq1}), which implies (\ref{eq2}) after taking an essential supremum with respect to $\mathcal H^s$ on the left hand side. 

\noindent Since $\Omega$ is open, there exists a closed ball $B_R(x)$ contained within $\Omega$.

\noindent We make use of the earlier mentioned spherical coordinate system - without loss of generality, after a rotation and a translation, we may assume $x = 0$, and $\nabla f(0) = (0, \dots, 0, |\nabla f(0)|)$. 

\noindent Now let $G \subset \Omega$ be the set of points of differentiability of $f$, and set $E = \Omega \setminus G$. We will associate to each $\varepsilon > 0$ a sequence of nested hypersurfaces $B_R (0) = H^\varepsilon_0 \supset H^\varepsilon_1 \supset \dots \supset H^\varepsilon_{n-1}$; and points $x\varepsilon_i \in H\varepsilon_i$ such that the following five conditions hold: 

\noindent i) Each $H^\varepsilon_i$ is the intersection of an $n-i$ dimensional linear subspace through the origin with $B_R(0)$.

\noindent ii) The intersection $H^\varepsilon_i \cap E$ is $\mathcal H^{n-1-i}$-null.

\noindent iii) For every $1 \leq i \leq n-1$, we have $\|\nabla f\|_{L^\infty (\mathcal H^{n-i}, H^\varepsilon_i)} \leq \|\nabla f\|_{L^\infty (\mathcal H^{n-i+1}, H^\varepsilon_{i-1})}$.

\noindent iv) For every $1 \leq i \leq n-1$, we have $|x^\varepsilon_i -x^\varepsilon_{i-1}| < \varepsilon$.

\noindent v) We have $x^\varepsilon_0 =R\frac{\nabla f(0)}{|\nabla f(0)|}$.

\noindent We proceed inductively. For the base case, we take $H^\varepsilon_0 := B_R (0)$, and $x^\varepsilon_0 := R\frac{\nabla f(0)}{|\nabla f(0)|}$ as specified by condition (v). Note that conditions (i) to (iii) are satisfied by these choices as well.

\noindent For the inductive step, let $1 \leq i \leq n-1$ be an integer and assume $H^\varepsilon_0, \dots, H^\varepsilon_{i-1}$; and $x^\varepsilon_0, \dots, x^\varepsilon_{i-1}$ satisfying conditions (i) to (v) have already been defined. Since $H^\varepsilon_{i-1}$ is a hyperplane isometric to the ball of radius $R$ around $0$ in $\mathbb R^{n-i+1}$, we may apply our spherical coordinate system, taking $H_{i-1}$ itself as the ambient space for the coordinates. After a rotation of the coordinate system we may assume that $H^\varepsilon_{i-1}$ is the closed ball of radius $R$ around $0$ in $\mathbb R^{n-i+1}$, and that $x^\varepsilon_{i-1}$ is the point $(0, \dots, 0, |x^\varepsilon_{i-1}|)$ in Euclidean coordinates, so that $x^\varepsilon_{i-1} = \Phi(|x^\varepsilon_{i-1}|, \frac{\pi}{2}, \dots, \frac{\pi}{2}, \frac{3\pi}{2})$.

\noindent Denote as usual $S_z$ the hypersurface obtained by fixing the last angular coordinate at $z$. By Lemma \ref{slicing} and Lemma \ref{iterated}, for every $r > 0$, $\mathcal H^1$ almost every choice of $z$ in $(\frac{3\pi}{2} - r, \frac{3\pi}{2}+r)$ satisfies conditions (ii) and (iii) above, and by Lemma \ref{hypersurface}, every choice satisfies (i). 

\noindent We fix some $r > 0$, whose precise value will be determined later. Fix some $z \in (\frac{3\pi}{2} - r, \frac{3\pi}{2} + r)$ satisfying (ii) and (iii). We then set $H^\varepsilon_i = S_{z} \cap B_R(0)$, and set $x^\varepsilon_i$ to be the orthogonal projection of $x^\varepsilon_{i-1}$ onto $H^\varepsilon_i$.

\noindent By construction, properties (i) to (iii) are satisfied by $H^\varepsilon_i$ and $x^\varepsilon_i$. It remains to check property (iv) for $x^\varepsilon_i$. But by standard properties of orthogonal projection, we have
$$\begin{aligned}|x^\varepsilon_i - x^\varepsilon_{i-1}| &= \inf_{h \in H^\varepsilon_i} |h^\varepsilon - x^\varepsilon_{i-1}| \\
& \leq\left |\Phi(|x^\varepsilon_{i-1}|, \frac{\pi}{2}, \dots, \frac{\pi}{2}, z) - \Phi(|x^\varepsilon_{i-1}|, \frac{\pi}{2}, \dots, \frac{\pi}{2}, \frac{3\pi}{2}) \right|\\
& \leq \text {Lip}(\Phi, B_R (0) \cap H_i) |z - \frac{3\pi}{2}| \\
& \leq \text{Lip}(\Phi, B_R (0) \cap H_i) r \\
& \leq \varepsilon,\end{aligned}$$
once $r$ is chosen less than $\frac{\varepsilon}{\text{Lip}(\Phi, B_R (0) \cap H_i)}$.

\noindent We have thus completed our construction satisfying conditions (i) to (v). Now given $\varepsilon > 0$, set $y_{\varepsilon} = \frac{x^\varepsilon_{n-1}}{|x^\varepsilon_{n-1}|}$, where $x_{n-1}$ is defined as above. We claim that
\begin{equation}\label{limit}|\nabla f(0)| = \lim_{\varepsilon \to 0} \langle \nabla f(0), y_\varepsilon \rangle.\end{equation}

\noindent Indeed. let $0 < \delta < \frac{1}{2}$ be arbitrary. Fix $\varepsilon = \frac{R\delta}{n-1}$ and a set of hyperplanes $H^\varepsilon_i$ and points $x^\varepsilon_i \in H_i$ associated to $\varepsilon$ as defined earlier. We have, by the triangle inequality
$$|x^\varepsilon_0 - x^\varepsilon_{n-1}| \leq |x^\varepsilon_0 - x^\varepsilon_1| + |x^\varepsilon_1 - x^\varepsilon_2| + \dots |x^\varepsilon_{n-2} - x^\varepsilon_{n-1}| \leq (n-1)\varepsilon.$$

\noindent Since $x^\varepsilon_0 = R\frac{\nabla f(0)}{|\nabla f(0)|}$, by the triangle inequality we have
$$R - (n-1)\varepsilon \leq |x_{n-1}| \leq R + (n-1)\varepsilon,$$
and so
$$(1-\delta)\leq \frac{|x_{n-1}|}{R} \leq (1+\delta)$$
by the definition of $\varepsilon$.

\noindent Thus we estimate
$$\begin{aligned}\left |\frac{\nabla f(0)}{|\nabla f(0)|} - \frac{x_{n-1}}{|x_{n-1}|} \right| & = \frac{1}{R}\left |R\left (\frac{\nabla f(0)}{|\nabla f(0)|} \right )- \frac{x_{n-1}}{|x_{n-1}|/R} \right |\\
&= \frac{1}{R} \max \left (|x_0 - \frac{x_{n-1}}{1-\delta}|, |x_0 -\frac{x_{n-1}}{1+\delta}| \right )\\
&\leq\frac{1}{R}\left ( |x_0 - x_{n-1}| + \frac{\delta|x_{n-1}|}{1-\delta}\right )\\
&\leq \frac{1}{R}\left (R\delta + \frac{\delta}{1-\delta}(1+\delta) R \right ) \\
&\leq \delta + 2\delta(1+\delta)\\
&\leq 4\delta.\end{aligned}.$$
From the above estimate, we deduce
$$\begin{aligned}\left ||\nabla f(0)| - \langle\nabla f(0), \frac{x_{n-1}}{|x_{n-1}|} \rangle \right | 
&= \left |\langle \nabla f(0), \frac{\nabla f(0)}{|\nabla f(0)|} \rangle -\langle \nabla f(0), \frac{x_{n-1}}{|x_{n-1}|} \rangle \right |\\
&=\left |\langle \nabla f(0), \frac{\nabla f(0)}{|\nabla f(0)|}- \frac{x_{n-1}}{|x_{n-1}|} \rangle \right |\\
&\leq 4\delta|\nabla f(0)|, \end{aligned}$$

\noindent by Cauchy-Schwarz. This tends to $0$ as $\delta \to 0$, and hence also $\varepsilon \to 0$. Thus \ref{limit} is proved.

\noindent We recall that by construction, $H\varepsilon_{n-1}^\varepsilon \cap E$ has $\mathcal H^0$ measure $0$, i.e. $f$ is differentiable everywhere over $H_{n-1}^\varepsilon$. Further, we have the isometry $ H^\varepsilon_{n-1} \cong \mathbb R \cap B_R (0)$, so by Proposition \ref{oned} applied to the restriction of $f$ to $H^\varepsilon_{n-1}$, we have
\begin{equation}\label{argh}\sup_{H^\varepsilon_{n-1}} |\nabla_{H^\varepsilon_{n-1}} f| \leq \|\nabla_{H^\varepsilon_v} f\|_{L^\infty(\mathcal H^{1}, H^\varepsilon_{n-1})},\end{equation} 
where $\nabla_{H_{n-1}^\varepsilon} f$ denotes the gradient of $f$ over the hypersurface $H_k^{\varepsilon}$. Denoting by $\mathcal B$ an orthogonal basis of $H^{\varepsilon}_k$, we recall that the norm of the gradient $|\nabla_{H^\varepsilon_{n-1}}f(0)|$ can be expressed as

$$\sup_{v \in \text{Span}(\mathcal B), |v| = 1} \langle \nabla f(0), v\rangle.$$

\noindent In particular, since $x^\varepsilon_{n-1} \in H^\varepsilon_{n-1}$ we have $y_\varepsilon \in \text{Span}(\mathcal B)$, and so 

\begin{equation}\label{argh2}\langle \nabla f(0), y_\varepsilon \rangle \leq |\nabla_{H^\varepsilon_{n-1}}f(0)|\end{equation}

\noindent We note that as a general fact, at any point $y$ of differentiability, we have $|\nabla_{H^{\varepsilon}_k} f(y)| \leq |\nabla f(y)|$, as the former gradient is taken over a strictly smaller set of directions. Thus
\begin{equation}\label{argh3}\|\nabla_{H^\varepsilon_{n-1}} f\|_{L^\infty(\mathcal H^{1}, H^\varepsilon_{n-1})} \leq \|\nabla f\|_{L^\infty(\mathcal H^{1}, H^\varepsilon_{n-1})} \end{equation}
\noindent Repeatedly applying condition (iii) above, we have
\begin{equation}\label{argh4}\begin{aligned} \|\nabla f\|_{L^\infty(\mathcal H^{1}, H^\varepsilon_{n-1})} &\leq \|\nabla f\|_{L^\infty(\mathcal H^{2}, H^\varepsilon_{n-2})} \\
&\ \,\vdots \\
&\leq \|\nabla f\|_{L^\infty(\mathcal H^{n}, H^\varepsilon_0)}\\
& \leq \|\nabla f\|_{L^\infty (\mathcal H^n, \Omega)}. \end{aligned}\end{equation}
\noindent We finally conclude
$$\begin{aligned}|\nabla f(0)| & = \lim_{\varepsilon \to 0_+} \langle \nabla f(0), y_\varepsilon \rangle \\
&\leq \lim_{\varepsilon \to 0_+} |\nabla_{H^\varepsilon_{n-1}} f(0)| \\
&\leq \lim_{\varepsilon \to 0_+} \|\nabla_{H^\varepsilon_{n-1}} f\|_{L^\infty(\mathcal H^{1}, H^\varepsilon_{n-1})} \\
&\leq \lim_{\varepsilon \to 0_+} \|\nabla f\|_{L^\infty(\mathcal H^{1}, H^\varepsilon_{n-1})}\\
& \leq  \|\nabla f\|_{L^\infty (\mathcal H^n, \Omega)}, \end{aligned}$$

where we have applied in sequence (\ref{limit}), (\ref{argh2}), (\ref{argh}), (\ref{argh3}), then (\ref{argh4}). \end{proof}
\section{Proof of remaining claims}
In this section, we present the proofs of Proposition $\ref{corollary2}$ and Corollary \ref{corollary3}. 

\begin{proof}[Proof of Proposition \ref{corollary2}] 

Since $f_n \to f$ in $W^{1, \infty}$, we have that 
$$\limsup_{n, m \to \infty} \, \underset{\Omega}{\esssup} \, |\nabla f_n - \nabla f_m| \to 0,$$
and thus, by Theorem \ref{maintheorem}, also
$$\limsup_{n, m \to \infty} \,  \sup_{\Omega} |\nabla f_n - \nabla f_m| \to 0.$$

\noindent That is, the sequence $\{\nabla f_n\}$ is Cauchy in supremum norm and hence converges uniformly to some $g: \Omega \to \mathbb R^n$. We claim now that $f$ is differentiable with $\nabla f = g$.

First, we show that the functions $f_n$ are uniformly differentiable - that is, for each $x \in \mathbb R^n$, we have that for every $\varepsilon > 0$ there exists $\delta > 0$ such that for all $n \in \mathbb N$,
\begin{equation}\label{unifdif}\frac{|f_n (x+y) - f_n (x) - \nabla f_n(x) \cdot y|}{y} < \varepsilon\end{equation}
whenever $|y| < \delta$. To this end, let $\varepsilon > 0$ be arbitrary. We pick $N \in \mathbb N$ so large such that
$$\sup_{n, m \geq N} \sup_{\Omega} |\nabla f_n - \nabla f_m| \leq \varepsilon.$$
Then for all $n \geq N$, we have that
$$\begin{aligned}&|f_n (x+y) - f_n (x) - \nabla f_n(x) \cdot y| \\
& \leq |f_n (x+y) - f_n (x) - (f_N (x+y) - f_N (x))| + |f_N (x+y) - f_N (x) - \nabla f_N (x) \cdot y| \\
&+ |(\nabla f_N(x) - \nabla f_n (x)) \cdot y| \\
& \leq \sup_\Omega |\nabla(f_n - f_N)(z)| \,  |y| + |f_N (x+y) - f_N (x) - \nabla f_N (x) \cdot y| + \varepsilon |y| \\
& \leq 2 \varepsilon |y| + |f_N (x+y) - f_N (x) - \nabla f_N (x) \cdot y|,\end{aligned}$$
so
$$\frac{|f_n (x+y) - f_n (x) - \nabla f_n (x)\cdot y|}{|y|} \leq 2\varepsilon + \frac{|f_N (x+y) - f_N (x) - \nabla f_N (x)|}{|y|}.$$
By differentiability of $f_N$, there is some $\delta_N$ such that the last term abovw is less than $\varepsilon$ for all $|y| < \delta_N$. Similarly, for each $0 \leq i < N$, there exists some $\delta_i$ such that 
$$\frac{|f_i (x+y) - f_i (x) - \nabla f_i (x) \cdot y|}{|y|} < \varepsilon$$
for all $|y| < \delta_i.$ Since $\varepsilon > 0$ was arbitrary, setting $\delta := \min_{0 \leq i \leq N} \delta_i$, we conclude.

Now we show the differentiability of $f$. Let $x \in \Omega$, $\varepsilon > 0$ be arbitrary. By the earlier discussion, there exists some $\delta > 0$ such that

$$\frac{|f_n (x+y) - f_n (x) - \nabla f_n (x) \cdot y|}{|y|} \leq \varepsilon$$
for all $y < \delta$.

Fix a small parameter $0 < h < 1$ which we will later send to $0$. Fix $n \in \mathbb N$ so large such that $\sup_{\Omega} |\nabla f_n - \nabla f| < \varepsilon$ and $\sup_{\Omega} |f_n - f| < h \delta \varepsilon$. Then for all $y \in B_{\delta}(x) \setminus B_{h\delta}(x)$, we compute using (\ref{unifdif}):
$$\begin{aligned}&|f(x+y) - f(x) - g(x) \cdot y|\\&\leq |f(x+y) - f_n (x+y)| + |f(x) - f_n (x)| + |f_n (x+y) - f_n (x) - \nabla f_n (x) \cdot y| \\
&+ |(\nabla f_n (x) - g) \cdot y|,\end{aligned} $$
so
$$\begin{aligned}\frac{|f(x+y) - f(x) - g(x) \cdot y|}{|y|} & \leq 2 \varepsilon + \varepsilon + \varepsilon \\
& = 4\varepsilon.\end{aligned}$$
We reiterate that the above holds for all $y \in B_{\delta} \setminus B_{h\delta}$.

Now fix a sequence of positive numbers $h_i \to 0_+$. Repeating the above argument with $h_i$ in place of $h$, and consequently $n_i$ in place of $n$, we obtain that for all $i$, and for all $y \in B_\delta (x) \setminus B_{h_i \delta}(x)$, we have
$$\frac{|f(x+y) - f(x) - g(x) \cdot y|}{|y|} < 4\varepsilon.$$
Since $h_i \to 0$, the above actually holds for all $y \in B_\delta(x) \setminus \{0\}$, which, since $\varepsilon > 0$ was arbitrary, shows the differentiability of $f$ at $x$ with derivative $g(x)$, as was to be shown.
\end{proof}

\begin{proof}[Proof of Corollary 3]
From the second part of Proposition $2$, we have that $\nabla f_n$ is Cauchy in uniform norm, hence $\nabla f_n \to \nabla f$ uniformly. As uniform limits of continuous functions are continuous, it follows that $\nabla f$ is continuous, as claimed.
\end{proof}
\section{Some follow-up questions}
We conclude with two natural follow-up questions. First, can we extend the main Theorem \ref{maintheorem} to sets of differentiability that have Hausdorff dimension $s > n-1$?

\noindent Second, we note that we did not use the full strength of Theorem \ref{maintheorem} in Corollary \ref{corollary2}. Can we strengthen the corollary to show that $W^{1, \infty}$ convergence preserves differentiability almost everywhere with respect to $s$-dimensional Hausdorff measure, for $0 \leq s \leq n-1$? The present corollary on everywhere differentiability corresponds to the case $s = 0$.

\printbibliography

 \end{document}